\documentclass[12pt]{article}
\usepackage{amssymb}
\usepackage[dvipdfm]{graphicx}
\usepackage{amsfonts}
\usepackage{latexsym}
\usepackage{amsthm}
\usepackage{amsmath}

\usepackage{amssymb, amscd}

\usepackage{url}

\newtheorem{problem}{Problem}[section]
\newtheorem{theo}[problem]{Theorem}
\newtheorem{rem}[problem]{Remark}
\newtheorem{prob}[problem]{Problem}
\newtheorem{defin}[problem]{Definition}
\newtheorem{prop}[problem]{Proposition}

\newtheorem{exam}[problem]{Example}

\newtheorem{conj}[problem]{Conjecture}

\begin{document}

 \title{A glimpse into continuous combinatorics\\ of posets, polytopes, and matroids}

\author{Rade  T.\ \v Zivaljevi\' c\\ {\small Mathematical Institute}\\[-2mm] {\small SASA, Belgrade}\\[-2mm]
 {\small rade$@$mi.sanu.ac.rs} }
\date{November 4, 2015}

\maketitle 

\begin{abstract}
Following \cite{Zi98} we advocate a systematic study of continuous
analogues of finite partially ordered sets, convex polytopes,
oriented matroids, arrangements of subspaces, finite simplicial
complexes, and other combinatorial structures. Among the
illustrative examples reviewed in this paper are an Euler formula
for a class of `continuous convex polytopes' (conjectured by Kalai
and Wigderson), a duality result for a class of `continuous
matroids', a calculation of the Euler characteristic of ideals in
the Grassmannian poset (related to a problem of Gian-Carlo Rota),
an exposition of the `homotopy complementation formula' for
topological posets and its relation to the results of Kallel and
Karoui about `weighted barycenter spaces' and a conjecture of
Vassiliev about simplicial resolutions of singularities. We also
include an extension of the index inequality (Sarkaria's
inequality) based on interpreting diagrams of spaces as continuous
posets.
\end{abstract}

\section{Introduction}

The idea of blending continuous and discrete mathematics into a
single `ConCrete' mathematics is far from being surprising or new.
Moreover, there seem to exist many different ways to carry on this
project, see for example \cite{GKP} (where calculus and
combinatorics interact in a fascinating way) and \cite{KR97}
(where the analogies between invariant measures on polyconvex sets
and measures on order ideals of finite partially ordered sets are
investigated). These are not isolated examples as exemplified by
papers \cite{AD12}, \cite{KW08}, \cite{Zi98}, which all address
some aspect of the problem of studying continuous objects from
discrete point of view or vice versa.

\medskip
Following into footsteps of \cite{Zi98}, in this paper are
collected some of the authors unpublished observations (and
impressions) about topological aspects of the problem of blending
discrete and continuous mathematics.

\medskip
In Section~\ref{sec:cont-polytopes} we explore (following Kalai
and Wigderson \cite{KW08}) the idea of studying convex bodies as
`continuous convex polytopes' (with continuous families of faces,
`continuous $f$-vector', etc.). The central result is an
Euler-style formula (Theorem~\ref{thm:Euler}) established for a
class of `tame convex bodies'.

\medskip
Section~\ref{sec:matroids} offers a brief treatment of `continuous
matroids'.  The central observation
(Proposition~\ref{prop:duality}) is that a simple convexity
argument can be used to show that continuous matroids, as
introduced in Section~\ref{sec:matroids}, have naturally defined
dual matroids satisfying a version of matroid duality.

\medskip
Topological partially ordered sets (or continuous posets for
short) are the most developed and possibly the most useful class
of `continuous-discrete' objects analyzed in this paper. In
Section~\ref{sec:GCRota} we focus on the Grassmannian topological
poset and show (Theorem~\ref{thm:Rota-answer}) its connection with
one of the problems of Gian-Carlo Rota from \cite{Rota98}. The
role of topological posets in the far reaching theory of
resolution of singularities (as founded and developed by Victor
Vassiliev \cite{Vas97}) is illustrated in
Section~\ref{sec:topo-poset}. Following \cite{Zi98} here we give a
brief exposition of the `homotopy complementation formula' for
topological posets. Among central examples is the configuration
poset ${\rm exp}_n(X)$ and one of the highlights is an exposition
of its relation to the `barycenter spaces' \cite{KK11} of Kallel
and Karoui and its connection to a conjecture of Vassiliev
(proposed on the conference `Geometric Combinatorics', MSRI
Berkeley, February 1997).

\medskip
Diagrams of spaces and their homotopy colimits appear in
Section~\ref{sec:sarkaria}. Here we illustrate how the
`continuous-discrete' point of view naturally leads to a useful
extension of the index inequality (Sarkaria's inequality from
\cite{Ziv} and \cite{Ma03}) to the case of diagrams of spaces
(Proposition~\ref{prop:Sarkaria-Zivaljevic}).

\section{Continuous polytopes}\label{sec:cont-polytopes}

Each convex body $K\subset \mathbb{R}^d$ can be interpreted as a
`continuous polytope'  (or $C$-polytope for short) with (possibly)
non-discrete families $F_k(K)$ of its $k$-dimensional faces. By
definition $A\in F_k(K)$ is a $k$-dimensional face of $K$ if $A$
is a $k$-dimensional closed convex set and,
\begin{enumerate}
\item[$\bullet$] for each line segment $[a,b]\subset K$ if
$(a,b)\cap A\neq\emptyset$ then $\{a,b\}\subset A$.
\end{enumerate}

It easily follows from the definition that if $A$ is a face of $B$
and $B$ is a face of $C$ then $A$ is a face of $C$. The set
$F_k(K)$ of all $k$-dimensional faces is naturally topologized by
the Hausdorff metric on the set of closed subsets of
$\mathbb{R}^d$.
\begin{defin}\label{def:face-space}
The disjoint union $\mathcal{F}(K) = \coprod_{k=0}^d F_k(K)$ is
referred to as the face-space of the convex body (continuous
polytope) $K$. The associated topological face poset is
$\mathcal{F}_K = (\mathcal{F}(K), \prec)$ where $A\prec B$ is the
containment relation $A\subseteq B$.
\end{defin}
A face $A\in \mathcal{F}(K)$ is `exposed' if $A = K\cap H$ for
some supporting hyperplane $H$ of $K$. Let $F_k^{exp}(K)\subset
F_k(K)$ be the space of all $k$-dimensional exposed faces of $K$
and $\mathcal{F}^{exp}(K)$ the associated space of all exposed
faces of $K$.

\medskip
If $A\in \mathcal{F}^{exp}(B)$ and  $B\in \mathcal{F}^{exp}(C)$
then it is not necessarily true that $A\in \mathcal{F}^{exp}(C)$.
For example an extremal point $a$ of $K$ which is not exposed such
that $[a,b]\in \mathcal{F}_1^{exp}(K)$ for some $b$ is an example
of a $0$-dimensional face with this property. Note that this is
not an isolated phenomenon since the Minkowski sum $K = O + P$ of
a smooth convex body $O$ and a convex polytope $P$ always have
points of this type.

\begin{figure}[hbt]
\centering
\includegraphics[scale=0.6]{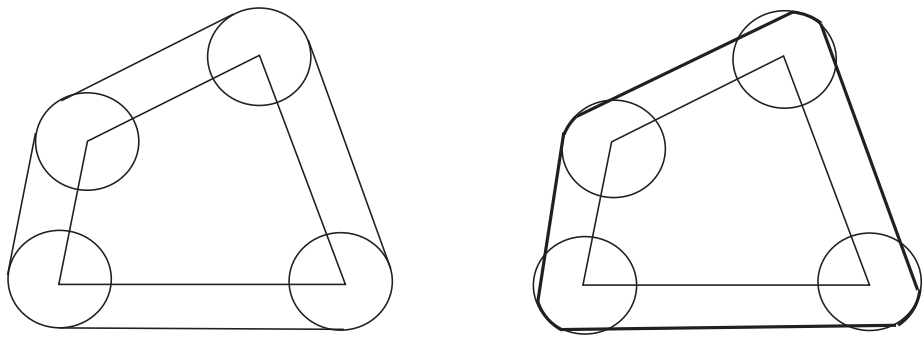}
\caption{Regularized Minkowski sum.} \label{figa:minko}
\end{figure}

\begin{enumerate}
\item[$\bullet$] The fact that $\mathcal{F}(K)$ is apparently
better behaved (as a topological poset) then the space
$\mathcal{F}^{exp}(K)$ is the reason why  we work mainly with
$\mathcal{F}(K)$.
\end{enumerate}
Let us make an empirical observation (without a formal proof) that
the Minkowski sum $K=O+P$ can modified (truncated, regularized) to
a convex body $K'$ which has better behaved facial structure and
which is often topologically similar to the original body $K$ in
the sense that $F_k(K)$ and $F_k(K')$ have the same homeomorphism
type, (Figure~\ref{figa:minko}).

\begin{prob}{\rm It would be useful to have a theorem providing a
regularization result illustrated in Figure~\ref{figa:minko} for
as large class of compact convex bodies as possible. More
precisely the problem is to construct, for a given compact convex
body $P$, a new convex body $Q$ such that,
\begin{enumerate}
 \item[(a)] each face of $Q$ is exposed,  $\mathcal{F}^{exp}(Q) = \mathcal{F}(Q)$;
 \item[(b)] $F_k(P)$ and $F_k(Q)$ are homeomorphic
 (homotopic) for each $k=0,\ldots, d$.
\end{enumerate}
 }
\end{prob}

\subsection{Tame continuous polytopes}
\label{sec:tame}

Our main objective in Section~\ref{sec:cont-polytopes} is
Theorem~\ref{thm:Euler} which confirms the Kalai-Wigderson
conjecture (Conjecture~\ref{conj:KW}) in the class of `tame
continuous polytopes'. Recall that the {\em Steiner centroid} is
the continuous selection $SC: \mathcal{K}_d\rightarrow
\mathbb{R}^d$ of a point from each compact convex set $A\in
\mathcal{K}_d$ which is Minkowski additive and invariant with
respect to Euclidean motions \cite{S71} (see also \cite{Zi89} for
some related facts and observations).

\begin{defin}\label{def:face-regular}
We say that a convex body $K\subset \mathbb{R}^d$ with compact
face-space $\mathcal{F}(K)$ {\rm
(Definition~\ref{def:face-space})} is $k$-face regular or $k$-face
tame if,
\begin{enumerate}
 \item[{\em (1)}] The collection $\{E_A\}_{A \in F_k(K)}$ of
 $k$-dimensional `tangent spaces' of $K$ at the $k$-dimensional faces is a vector
 bundle $\pi_k : \mathcal{E}_k \rightarrow F_k(K)$  over $F_k(K)$;
 \item[{\em (2)}] Let $\mathcal{C}_k = \bigcup\{\mbox{\rm relint}(A) \mid A\in
 F_k(K)\}$ be the union of relative interiors of all
 $k$-dimensional faces of $K$ and $\widehat{\mathcal{C}}_k$ its
 one-point compactification. Then the space
 $\widehat{\mathcal{C}}_k$ and the Thom space ${\rm Thom}(\mathcal{E}_k)$ of the bundle
 $\mathcal{E}_k$ are homeomorphic.
 \end{enumerate}
A convex body $K$ is {\em `face lattice tame'} or simply tame if
it is $k$-face regular for each $0\leq k\leq \mbox{\rm dim}(K)$.
\end{defin}

The conditions (1) and (2) in Definition~\ref{def:face-regular}
may require a little clarification. For $A\in F_k(A)$ the affine
span $\mbox{\em aff}(A)$ of $A$ is naturally a vector space with
$0 = 0_A$ as the origin; more explicitly $E_A$ is the vector
subspace of $\mathbb{R}^d$ obtained by translating $\mbox{\em
aff}(A)$ by the vector $-SC(A)$. The condition (1) says that this
family of vector spaces is locally trivial which means that
$\mathcal{E}_k := \bigcup_{A\in F_k(K)}~\{A\}\times E_A\subset
F_k(K)\times \mathbb{R}^d$ is a total space of a genuine vector
bundle over $F_k(K)$.

\medskip
The condition (2) says that we are allowed to treat individual,
$k$-dimensional, closed convex sets $A\in F_k(K)$ as `discs' in
$E_A$ and (more importantly) the union of relative interiors of
all $A\in F_k(K)$ as the total space of the open disc bundle
associated to the bundle $\mathcal{E}_k$.

\begin{prob}{\rm
It would be certainly nice to have a description of general
classes of convex bodies which are `face lattice tame' in the
sense of Definition~\ref{def:face-regular}. } \label{prob:tame}
\end{prob}

\begin{exam}\label{exam:wild} {\rm
In the direction opposite to Problem~\ref{prob:tame} one can
search for the simplest examples of $C$-polytopes which are `wild'
in the sense that they violate either (1) or (2) in
Definition~\ref{def:face-regular}. A $3$-dimensional example
arises by taking the convex hull ${\rm conv}(D\cup I)$ where $D =
\{(x,y)\mid (x-1)^2+y^2\leq 1\}$ is a unit disc in the $xy$-plane
and $I$ is a vertical segment on the $z$-axis which contains the
origin in its interior. }
\end{exam}

\subsection{Euler formula for continuous polytopes}\label{sec:Euler}

Kalai and Wigderson conjectured in \cite[Conjecture~6]{KW08} the
following Euler type formula for continuous polytopes. Here and
elsewhere $\chi(X)$ is the Euler characteristic of the space $X$.

\begin{conj}\label{conj:KW} Suppose that $K$ is a convex body in $\mathbb{R}^d$
and let $F_k(K)$ be the space of all $k$-dimensional faces of $K$
with the topology induced by the Hausdorff metric. Assume that
$F_k(K)$ is compact. Then,
\begin{equation}\label{eqn:Kalai-Wigderson}
\sum_{k=0}^{d-1} (-1)^k \chi(F_k(K)) = \chi(S^{d-1}) = 1 +
(-1)^{d-1}.
\end{equation}
\end{conj}

\begin{theo}\label{thm:Euler}
Suppose that $K$ is a convex body which is {\em face lattice tame}
in the sense of Definition~\ref{def:face-regular}. Then,
\begin{equation}\label{eqn:Kalai-Wigderson-2}
\sum_{k=0}^{d-1} (-1)^k \chi(F_k(K)) = \chi(S^{d-1}) = 1 +
(-1)^{d-1}.
\end{equation}
In other words Conjecture~\ref{conj:KW} is true if for each $k$
the space $F_k(K)$ is essentially the base space of a naturally
associated vector bundle $\mathcal{E}_k$
(Definition~\ref{def:face-regular}).
\end{theo}

\medskip\noindent
{\bf Proof:} Let $$\mathcal{F}_{k}= \mathcal{F}_{k}(K) :=
\bigcup_{j=0}^k F_j(K)$$ be the union of all $j$-dimensional faces
of $K$ for $j=0,\ldots, k$.

\medskip
By definition
 $$
\mathcal{F}_{k}\setminus \mathcal{F}_{k-1} = \bigcup\{\mbox{\rm
relint}(A) \mid A\in F_k(K)\}
$$
is the union of relative interiors of all $k$-dimensional faces of
$K$ and there is commutative square,
\begin{equation}
\begin{CD}
\mathcal{F}_{k}\setminus \mathcal{F}_{k-1}  @>\alpha>>
\mathcal{E}_k\\
@V\pi_kVV   @VV\pi_kV \\
F_k(K) @>>=>  F_k(K)
\end{CD}
\end{equation}
where $\alpha$ is an inclusion map. By the tameness assumption
(Definition~\ref{def:face-regular}) the one-point compactification
$\widehat{\mathcal{C}}_k$ of $\mathcal{F}_{k}\setminus
\mathcal{F}_{k-1}$ is homeomorphic to the Thom space $T_k = {\rm
Thom}(\mathcal{E}_k)$ of the bundle $\mathcal{E}_k$.

\medskip
Let $\widetilde{\chi}(Y)$ be the reduced Euler characteristics of
a pointed space $Y$. By the Thom isomorphism theorem we know that
$\chi(F_k(K)) = (-1)^k\widetilde{\chi}(T_k)$ (here we took into
account the fact that the isomorphism shifts the dimension by
$k$). From the exact sequence of the pair $(\mathcal{F}_{k},
\mathcal{F}_{k-1})$ we deduce that,

\begin{equation}\label{eqn:chi}
\chi(\mathcal{F}_k) = \chi(\mathcal{F}_{k-1})
+\widetilde{\chi}(T_k) = \chi(\mathcal{F}_{k-1}) +
(-1)^k{\chi}(F_k(K)).
\end{equation}
Note that for $k=0$ the relation (\ref{eqn:chi}) reduces to
$\chi(\mathcal{F}_0) = \chi(F_0(K))$. By adding the equalities
(\ref{eqn:chi}) for $k=0,\ldots, d-1$ we obtain,
\begin{equation}
\sum_{k=0}^{d-1} (-1)^k \chi(F_k(K)) = \chi(\mathcal{F}_{d-1})
\end{equation}
and the Euler relation (\ref{eqn:Kalai-Wigderson-2}) follows from
the fact that $\mathcal{F}_{d-1} = \partial(K)\cong S^{d-1}$.
\hfill$\square$

\section{Continuous matroids}

`Continuous matroids' is another class of continuous objects
motivated by their discrete counterparts. The exposition in this
section is based on the unpublished manuscript \cite{Zi09}. The
central is Proposition~\ref{prop:duality} which shows that
continuous matroids, as introduced in Section~\ref{sec:matroids},
have naturally defined dual matroids satisfying a version of
matroid duality, cf. \cite[Lecture 6]{Zie} for a classical
treatment of the case of oriented matroids. The reader is referred
to \cite{AD12} for  an up-to-date treatment of continuous matroids
from a parallel point of view.

\subsection{Complex and quaternionic matroids}
\label{sec:matroids}

Suppose that $\mathbb{K}$ is one of the classical (skew) fields
$\mathbb{R}, \mathbb{C}$ or $\mathbb{H}$. Let $S=S_{\mathbb{K}}
\cong S^{d(K)-1}$ be the unit sphere in $\mathbb{K}$ and let
$\mathbb{K}^n$ be the $n$-dimensional vector space (left module)
over $\mathbb{K}$.

\begin{defin}\label{def:cross-polytope} A $\mathbb{K}$-cross polytope
in $\mathbb{K}^n$ is the convex body $\lozenge_{\mathbb{K}}^n$
defined as the convex hull
$$
\lozenge_{\mathbb{K}}^n := {\rm conv}\bigcup_{i=1}^n S_i
$$
where $S_i:= \{z\in \mathbb{K}^n \mid \vert z_i\vert = 1 \mbox{
{\rm and for all} } j\neq i, z_j = 0 \}$ is the unit sphere in the
$i^{\rm th}$ coordinate line.
\end{defin}

We see $\lozenge_{\mathbb{K}}^n$ as an example of a ``continuous''
polytope ($C$-polytope in the sense of
Section~\ref{sec:cont-polytopes}). Recall that a $C$-polytope is
simply a convex body that exhibits properties of both smooth
convex bodies and convex polytopes. Other examples of C-polytopes
include the ``continuous cyclic polytope'' defined as the convex
hull of the curve $\Gamma_n = \{(z,z^2,\ldots,z^n) \mid \vert
z\vert =1)\}$, or more generally convex hulls of embedded
manifolds, \cite{KW08}. Even more familiar examples (already met
in Section~\ref{sec:cont-polytopes}) are Minkowski sums of smooth
convex bodies and convex polytopes, in particular convex bodies of
the form $C = A\times Q\subset \mathbb{R}^m\times \mathbb{R}^n$
are good motivating examples of $C$-polytopes where $A$ is a
(possibly smooth) convex body in $\mathbb{R}^m$ and $Q\subset
\mathbb{R}^n$ a convex polytope.

\medskip
Summarizing  a $C$-polytope is just an ordinary convex body  $K$
portrayed as a some kind of a ``continuous convex polytope''. A
characteristic property of a $C$-polytope $K$ is that its face
poset (Definition~\ref{def:face-space}) is a continuous posets in
the sense of \cite{Zi98} (see also our
Section~\ref{sec:cont-posets}).

\begin{defin}\label{def:matroids}
Suppose that $K\subset \mathbb{R}^n$ is a $C$-polytope such that
$0\in {\rm int}(K)$. Let $\mathcal{F}_K$ be the associated
face-poset (Definition~\ref{def:face-space}). Let $L \subset
\mathbb{R}^n$ be a linear subspace. Then the $K$-matroid
$\mathcal{M}_K(L)$ of $L$ is by definition $\mathcal{M}_K(L) =
\{A\in \mathcal{F}_K \mid {\rm relint}(A)\cap L\neq\emptyset\}$. A
$\lozenge_{\mathbb{K}}^n$-matroid of $L$, where
$\lozenge_{\mathbb{K}}^n$ is the $\mathbb{K}$-cross polytope
described in Definition~\ref{def:cross-polytope}, is referred to
as a $\mathbb{K}$-matroid and denoted by
$\mathcal{M}_\mathbb{K}(L)$.
\end{defin}

\begin{exam}\label{exam:cross}{\rm
Suppose that $K=\lozenge_{\mathbb{R}}^n = \lozenge^n$ is the
``ordinary'' cross-polytope. Then the face poset
$\mathcal{F}_{\lozenge^n}$ (with $\emptyset = \hat{0}$ as the
minimum element) is isomorphic to the poset ${\mathcal
Sgn}_n=(\{-1,0,+1\}^n,\leq)$ of all sign vectors from the usual
theory of oriented matroids. By definition the $K$-matroid
$\mathcal{M}_{\mathbb{K}}(L)$, associated to a subspace $L\subset
\mathbb{R}^n$ is a {\em realizable oriented matroid} from the
standard theory or oriented matroids. Indeed,
$\mathcal{M}_{\mathbb{K}}(L)$ is essentially the collection of all
sign vectors $sgn(v)\in \{-1,0,+1\}^n$ for all $v\in L$.}
\end{exam}

\subsection{Sign vectors}

As already indicated in Example~\ref{exam:cross}, faces of a
$C$-polytope $K\subset \mathbb{R}^n$ should be understood as
generalized sign-vectors. In particular the map
\begin{equation}\label{eqn:sign-vect}
\nu : \mathbb{R}^n \rightarrow \mathcal{F}_K,
\end{equation}
which associates to a vector $v\in \mathbb{R}^n$ its sign
$\nu(v)$, is defined as the unique face $F\in \mathcal{F}_K$ such
that the ray $\rho(v):= \{\lambda v \mid \lambda\geq 0\}$ and
${\rm relint}(F)$ have a non-empty intersection.

\medskip
An ultimate justification for this definition is the fact (see
\cite[Theorem~18.2.]{R70}) that the collection $\{{\rm relint}(A)
\mid A\in \mathcal{F}_K\}$ is a partition of the $C$-polytope $K$.
In particular each ray $\rho(v)$ intersects precisely one of the
sets ${\rm relint}(A)$ for $A\in\mathcal{F}_K$.

\medskip
In analogy with the case of usual oriented matroids we call
$\nu(v)$ a $K$-sign or $K$-sign vector of $v$, in particular the
set of all vectors which share the same $K$-sign vector $F\in
\mathcal{F}_K$ is the (relatively open) cone, ${\rm cone}({\rm
relint}(F))$. The family of cones
$$
\mathcal{F} = \{{\rm cone}({\rm relint}(F)) \mid F\in
\mathcal{F}_K\}
$$
is a ``continuous-discrete'' fan in $\mathbb{R}^n$. Clearly one
could have started from the beginning with a $C$-fan, instead of
the $C$-polytope. However, at this stage it appears to be more
natural to explore in some detail the motivating examples so we
focus on the case of convex bodies with a particular emphasis on
bodies $\lozenge_{\mathbb{K}}^n$.

\subsection{Orthogonality and duality}

Suppose that $X$ and $Y$ are two vector spaces (left moduli) over
$\mathbb{K}$ and let $\langle\cdot,\cdot\rangle : X\times
Y\rightarrow \mathbb{K}$ be a non-degenerate bilinear form which
allows us to talk about orthogonality of vectors and sets in $X$
and $Y$. One could start with $C$-bodies $A\subset X$ and
$B\subset Y$, each with the corresponding families of $A$-matroids
and $B$-matroids, and try to develop a natural concept of duality
between these classes.

\medskip
Again, we temporary sacrifice generality and focus to the main
case of the convex body $\lozenge_{\mathbb{K}}^n$. Our objective
is to introduce an orthogonality relation for the associated
signed vectors which should lead to the duality of
$\mathbb{K}$-matroids.

\medskip
Let $\langle x, y \rangle = x_1\bar{y}_1+\ldots +x_n\bar{y}_n$ be
the standard Hermitian form on $\mathbb{K}^n$.

\begin{defin}\label{def:orthogonality}
We say that two signed vectors $a,b\in \mathcal{F}_{\mathbb{K}} =
\mathcal{F}_{\lozenge_{\mathbb{K}}^n}$ are orthogonal $a\perp b$,
if there exist vectors $x,y\in \mathbb{K}^n$ such that $a=\nu(x),
b=\nu(y)$ and $\langle x,y\rangle = 0$. Given a subset
$\mathcal{M}\subset  \mathcal{F}_{\mathbb{K}}$ let
$$
\mathcal{M}^\perp := \{b\in  \mathcal{F}_{\mathbb{K}} \mid
(\forall a\in \mathcal{M}) \, a\perp b\}.
$$
\end{defin}

The following statement,  claiming the compatibility of the
operations of the geometric and matroid dual, is possibly an
encouraging sign and a good omen for the theory of continuous,
complex and quaternionic matroids. For simplicity the
$\mathbb{K}$-matroid $\mathcal{M}_\mathbb{K}(V)$ of a vector space
$V$ is denoted by $\mathcal{M}(V)$.

\begin{prop}\label{prop:duality}
\begin{equation}\label{eqn:duality}
\mathcal{M}(V^\perp) = \mathcal{M}(V)^\perp .
\end{equation}
\end{prop}

\medskip\noindent
{\bf Proof:}  If $b\in \mathcal{M}(V^\perp)$ then $b = \nu(y)$ for
some $y\in V^\perp$. Hence $y\perp x$ for each $x\in V$ and
$b\perp a$ for each $a\in \mathcal{M}(V)$, which implies that
$b\in \mathcal{M}(V)^\perp$ and completes the proof of the
inclusion $\mathcal{M}(V^\perp) \subset \mathcal{M}(V)^\perp$.

\medskip
Let us prove the opposite inclusion $\mathcal{M}(V^\perp) \supset
\mathcal{M}(V)^\perp$ by contraposition. Suppose $b\notin
\mathcal{M}(V^\perp)$. Then $b$ corresponds to a face $F_b$ of
$\lozenge^n_{\mathbb{K}}$ and for each $x\in {\rm relint}(F_b), \,
x\notin V^\perp$, that is
$$
{\rm relint}(F_b)\cap V^\perp = \emptyset.
$$
By the separation principle for convex sets there exists a vector
$u\in \mathbb{K}^n$ such that,

\begin{enumerate}
\item[(1)] ${\rm Re}\langle z,u\rangle > 0$ for each $z\in {\rm
relint}(F_b)$;

\item[(2)] ${\rm Re}\langle z,u\rangle = 0$ for each $z\in
V^\perp$.
\end{enumerate}
Since $V^\perp$ is a left $\mathbb{K}$-module, it follows from (2)
that ${\rm Re}\langle z,\alpha u\rangle = 0$ for each $\alpha\in
\mathbb{K}$, which immediately implies,
\begin{enumerate}
\item[$(2^\sharp)$] $\langle z,u\rangle = 0$ for each $z\in
V^\perp$.
\end{enumerate}
From here we deduce $u\in V$. Let $a=\nu(u)$. Then $a\in
\mathcal{M}(V)$ and in light of (1), $b\not\perp a$ which finally
implies $b\notin \mathcal{M}(V)^\perp$. \hfill $\square$

\section{Continuous posets}\label{sec:cont-posets}

Continuous posets \cite{Vas91, Vas99}, \cite{Zi98} are perhaps the
most useful and widely applicable examples of continuous analogues
of discrete structures. One of the main and most interesting
examples of topological posets is the `Grassmannian poset'. For
the `order complex construction' (or the `flag-join' construction)
and all other undefined concepts and related results the reader is
referred to \cite{Vas91, Vas99} and \cite{Zi98}.

\begin{defin}\label{def:Grassman-poset}
The Grassmannian poset ${\cal G}_n(\mathbb{R}) =
(G(\mathbb{R}^n),\subseteq)$, is the disjoint sum  $$
G(\mathbb{R}^n) := \coprod_{i=0}^{n} G_i(\mathbb{R}^n) $$ where
$G_i(\mathbb{R}^n)$ is the manifold of all $i$-dimensional linear
subspaces of $\mathbb{R}^n$. The order in this poset is by
inclusion, $U\leq V$ iff $U\subseteq V$. Denote the minimum and
the maximum element in this poset by ${\hat 0}$ and ${\hat 1}$
respectively and let $\rho : {\cal G}_n(\mathbb{R}) \rightarrow
\mathbb{N}, \, L\mapsto {\rm dim}(L)$ be the rank function. The
poset $\widetilde{\cal G}_n(\mathbb{R}) := {\cal
G}_n(\mathbb{R})\setminus \{\hat 0,\hat 1\}$ is called the
truncated Grassmannian poset. Let $I\subset\widetilde{\cal
G}_n(\mathbb{R})$ be a closed order ideal (initial subset). The
order complex $\Delta(I)$, see \cite{Vas91} and \cite{Zi98}, is
defined as the subspace of the join
$$ G_1(\mathbb{R})\ast G_2(\mathbb{R}) \ast \ldots \ast G_{n-1}(\mathbb{R})$$
spanned by all flags in $I$.
\end{defin}

\begin{rem}\label{rem:remarkable}{\rm
It is a remarkable fact (see \cite{Vas91, Vas99, Zi98}) that the
order complex (flag-join) of the truncated Grassmannian poset
$\widetilde{\cal G}_n(\mathbb{R})$ is a sphere of dimension
${{n\choose 2}+n-2}$,
\begin{equation}\label{eqn:Grass-poset}
\Delta(\widetilde{\cal G}_n(\mathbb{R})) \cong S^{{n\choose
2}+n-2}.
\end{equation}
As an immediate consequence we obtain that,
\begin{equation}\label{eqn:Grass-poset-bis}
\Delta(\widehat{\cal G}_n(\mathbb{R})) \cong D^{{n\choose 2}+n-1}
\end{equation}
is a disc of dimension ${{n\choose 2}+n-1}$ where
$\widehat{\mathcal{G}}_n(\mathbb{R}) =
\mathcal{G}_n(\mathbb{R})\setminus \{\widehat{0}\}$. }
\end{rem}

\begin{defin}\label{def:chi-Grassmann}
Let $I\subset\widetilde{\cal G}_n(\mathbb{R})$ be a closed order
ideal in the Grassmannian poset and let $I_k = I\cap
G_k(\mathbb{R}^n)$. Then, $$\chi(I) = (\chi_1(I),\chi_2(I),\ldots,
\chi_n(I))$$ is referred to as the $\chi$-vector of the ideal $I$
where $\chi_k(I) = \chi(I_k)$.
\end{defin}

\subsection{Grassmann posets and a problem of Gian-Carlo Rota}
\label{sec:GCRota}

\begin{defin} Let $P$  be a topological poset equipped with a rank
function $\rho : P\rightarrow \mathbb{N}$. A $P$-complex is by
definition an order ideal $I$ in $P$. Let $I_m$ be the set of all
elements in $I$ of rank $m\in \mathbb{N}$. The $\chi$-vector of
the $P$-complex $I$ is by definition
$$
\chi_P(I) := (\chi(I_0),\chi(I_1),\ldots ,\chi(I_m),\ldots )
$$
where $\chi(X)$ is the Euler characteristic of $X$. For example if
$P$ is a simplex then $I$ is a simplicial complex and $\chi_P(I)$
is the usual $f$-vector of $I$. In this case there is a well-known
relation
\begin{equation}
\label{Euler_1} \chi(\Delta(I)) = f_0 - f_1 + f_2 - \ldots
\end{equation}
\end{defin}

Gian-Carlo Rota delivered on a joint meeting of the American
Mathematical Society and Mexican Mathematical Society (Oaxaca,
Mexico, December 1997) a lecture with a charming, provocative and
(in retrospective) saddening title `Ten Mathematics Problems I
will never solve'\footnote{Gian-Carlo Rota passed away on April
18, 1999.}, see \cite{Rota98} for the published version.

\medskip
Among Rota's problem is the Problem 7 (on Intrinsic volumes of
families of subspaces) where he formulates (in our language) the
problem of developing the theory of (finitely additive)
$O(n,\mathbb{R})$-invariant measures defined on the class of
closed order ideals of the Grassmann poset ${\cal
G}_n(\mathbb{R})$.

\medskip
Gian-Carlo Rota was guided by an analogy with the (simple and
well-understood) theory of $S_n$-invariant measures on the class
of order ideals in the posets of all subsets of the set $[n] =
\{1,2,\ldots,n\}$. In this case order ideals are nothing but the
simplicial complexes (on $[n]$ as the set of vertices) and Rota
relates the well known formula,
\begin{equation}\label{eqn:formula}
\chi(K) = f_0(K) - f_1(K) + f_2(K) +\ldots +(-1)^n f_n(K)
\end{equation}
to the fact that the  Euler characteristics $\chi$ is the unique,
$S_n$-invariant, finitely additive measure defined on simplicial
complexes.

\medskip
Rota concluded his description of Problem~7 by saying that `At
present, we cannot even get the Euler characteristics', in other
words  Rota pointed to the following special case of his
Problem~7,

\begin{prob}\label{prob:Rota}{\em
Find an analogue of  (\ref{eqn:formula}) for the class of closed
ideals in the poset ${\cal G}_n(\mathbb{R})$ of all linear
subspaces of a finite dimensional Hilbert space. }
\end{prob}

\medskip
The reader familiar with the results of Vassiliev
\cite{Vas91,Vas93} about the structure of the order complex of the
Grassmann posets (Remark~\ref{rem:remarkable}) will immediately
see that these results provide a key to the
Problem~\ref{prob:Rota}.

\medskip
The following theorem, from an unpublished manuscript \cite{Zi98b}
presented at the conference in Kotor-98,   provides an amusing
answer to the Problem~\ref{prob:Rota}.

\begin{theo}\label{thm:Rota-answer}
Let $I \subset \widetilde{\cal G}_n(\mathbb{R})$ be a closed ideal
in the truncated Grassmannian poset $\widetilde{\cal
G}_n(\mathbb{R})$ and let $\chi(I) = (\chi_1(I),\chi_2(I),\ldots)$
be the associated $\chi$-vector in the sense of
Definition~\ref{def:chi-Grassmann}. Then,
\begin{equation}
\label{Euler_3} \chi(\Delta(I)) = \chi_1(I) + \chi_2(I) -
\chi_3(I) - \chi_4(I) + \ldots + (\sqrt{-1})^{n^2+n+2}\chi_n(I)
\end{equation}
\end{theo}

\medskip\noindent
{\bf Proof:} The proof is similar to the proof of
Theorem~\ref{thm:Euler}. If $I_{\leq k} = I_1\cup I_2\cup \ldots
\cup I_k$ then $\Delta(I_{\leq k})\subset \Delta(I)$ and there is
an increasing filtration,
\begin{equation}\label{eqn:filtr-Grass-pos}
\Delta(I_{\leq 1})\subset \Delta(I_{\leq 2})\subset \ldots \subset
\Delta(I_{\leq k})\subset \ldots \Delta(I_{\leq n-1}) = \Delta(I)
\end{equation}
The central observation is that $\Delta(I_{\leq k})/
\Delta(I_{\leq k-1})\cong Thom(U_k)$ is the Thom space of a vector
bundle $U_k$ over $I_k$ of dimension $c_k = {k\choose 2} +k-1$.
Indeed, there is a set theoretic decomposition,
\begin{equation}\label{eqn:decomposition}
\Delta(I_{\leq k})\setminus \Delta(I_{\leq k-1}) = \biguplus_{V\in
I_k} \Delta(\widehat{\mathcal{G}}(V))\setminus
\Delta(\widetilde{\mathcal{G}}(V))
\end{equation}
where $\widehat{\mathcal{G}}(V)$ and $\widetilde{\mathcal{G}}(V)$
are respectively posets isomorphic to
$\widehat{\mathcal{G}}_k(\mathbb{R})$ and
$\widetilde{\mathcal{G}}_k(\mathbb{R})$ (described in
Remark~\ref{rem:remarkable}). These isomorphisms arise from
locally chosen isomorphisms  $V\cong \mathbb{R}^k$ (provided by
local trivializations of the canonical $k$-plane bundle over the
Grassmannian $G_k(\mathbb{R})$). In light of
Remark~\ref{rem:remarkable}
$\Delta(\widehat{\mathcal{G}}(V))\setminus
\Delta(\widetilde{\mathcal{G}}(V))$ is an open disc of dimension
$c_k = {k\choose 2} +k-1$. Moreover, upon closer inspection we see
that (\ref{eqn:decomposition}) is actually the total space of a
$c_k$-dimensional vector bundle associated to the canonical
$k$-plane bundle over $I_k$ induced from the canonical $k$-plane
bundle over $G_k(\mathbb{R}^n)$.

The proof is concluded in the same way as the proof of
Theorem~\ref{thm:Euler}. \hfill $\square$

\section{Topological posets}\label{sec:topo-poset}

\subsection{Poset resolution of $P$-singular spaces}

Vassiliev's ``Geometric resolutions of singularities''
\cite{Vas91, Vas92, Vas93, Vas97, Vas99} is a versatile and
powerful method for studying topology of {\em singular spaces} and
their complements. A substantial part of the theory can be
rephrased and fruitfully generalized in the language of
topological order complexes.

\medskip
A model example of a singular space is a subspace $X\subset
Fun(M,N)$ of some function space where $f\in X$ if and only if $f$
is degenerate in some (precisely defined) sense. Our objective is
to study the topology of the singular space $X$ by studying the
associated space $\widehat{X}$ obtained from $X$ by `resolving the
singularities'. The construction can be (somewhat informally)
summarized as follows.

\begin{enumerate}
 \item[$\bullet_1$] $X$ is a singular space, e.g.\ the space of
singular matrices, polynomials with multiple zeros, singular
knots, smooth functions that exhibit singularities of certain type
etc.
 \item[$\bullet_2$] There is a hierarchy of observed {\em
singularity types} which are naturally arranged in a topological
poset $(\mathcal{P}, \prec)$ where $p \prec p'$ means that the
singularity type  $p'$ is in some sense more complex than $p$.

 \item[$\bullet_2$] There is a map $\Phi : X  \rightarrow \mathcal{P}$
which associates to each point $x\in X$ its singularity type which
is `semi-continuous' in the sense that in the limit $x_n
\longrightarrow x$ the singularity type can only jump up in the
complexity (increase in $\mathcal{P}$).

 \item[$\bullet_4$] The $\mathcal{P}$-resolution of the singular space
 $X$ is the space,
$$
\widehat{X} := \bigcup_{x\in X}\{x\}\times
\Delta(\mathcal{P}_{\geq \Phi(x)})\subset X\times
\Delta(\mathcal{P}).
$$
It is expected that, as a consequence of semi-continuity of
$\Phi$, the space $\widehat{X}$ is a closed subset of $X\times
\Delta(\mathcal{P})$. Moreover we assume that the natural
projection $\pi : \widehat{X} \rightarrow X$ has contractible
fibers so (under mild assumptions) it is a homotopy equivalence.

 \item[$\bullet_5$] There is a global filtration of the poset
$\mathcal{P}$ (for example by a monotone rank function $\rho :
\mathcal{P} \rightarrow \mathbb{Z}$). This filtration induces a
filtration on $\widehat{X}$ which leads to a spectral sequence
computing the (co)homology of $X$.
\end{enumerate}

\medskip
The scheme described above appears to be so fundamental that the
very concept of a singular space may accordingly modified.  The
category of $\mathcal{P}$-singular spaces $Sing(\mathcal{P})$ is a
natural ambient for studying both the interesting
$\mathcal{P}$-singular spaces and the topological poset
$\mathcal{P}$ itself (where an object of $Sing(\mathcal{P})$ is
treated as some sort of a module (sheaf) over the ring (space)
$\mathcal{P}$).

\begin{defin}\label{def:P-sing-space}
Suppose that $P$ is a topological poset. A topological space $X$
is given a structure of a $P$-singular space if there is a map
$\Phi : X \rightarrow P$ which has (some or all of the) properties
$\bullet_1$ to $\bullet_5$. A morphism $X \dashrightarrow Y$
between two $\mathcal{P}$-singular spaces is a map over
$\mathcal{P}$ (commutative diagram) which preserves all the
associated structures listed in $\bullet_1$--$\bullet_5$. In
particular there is a map $\widehat{f} : \widehat{X}\rightarrow
\widehat{Y}$ of the associated $\mathcal{P}$-resolutions which
respects the filtrations described in $\bullet_5$ such that,
$$
\begin{CD}
\widehat{X} @>>>    \widehat{Y}\\
@V\pi VV     @VV\pi V   \\
 X   @>>>    Y
\end{CD}
$$
\end{defin}
The category $Sing(\mathcal{P})$ described in
Definition~\ref{def:P-sing-space} comes naturally with a functor
$S : Sing(\mathcal{P}) \longrightarrow SpecSeq$ from the category
of $\mathcal{P}$-singular spaces to the category of spectral
sequences. One should in principle be able to construct
simplifying test objects in $Sing(\mathcal{P})$ and use the
functor $S$ to detect (describe) particular (co)homology classes
(characteristic classes) of the $\mathcal{P}$-singular space under
consideration.

\subsection{Topological homotopy complementation formulas}
\label{sec:homotopy-complementation}

One of the central guiding principles of \cite{Zi98} is that
ideas coming from discrete combinatorics, properly interpreted and
generalized, can play a unifying and motivating role in the
analysis of topological (continuous) posets. The main example in
\cite{Zi98} of such a result about finite (discrete) posets is the
so called `Homotopy complementation formula' of Bj{\"o}rner and
Walker \cite{BW83} (HCF for short).

\begin{figure}[hbt]
\centering
\includegraphics[scale=0.7]{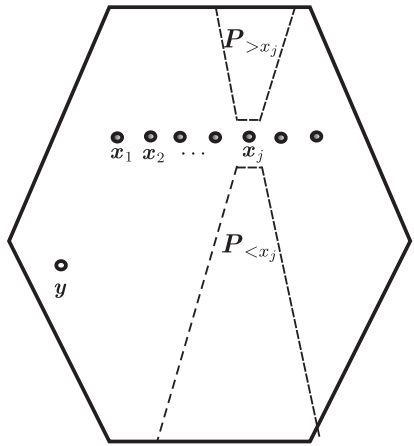}
\caption{Evaluation of the homotopy type of
$\Delta(P)/\Delta(P\setminus X)$.} \label{fig:poset}
\end{figure}
\noindent Suppose that $X = \{x_j\}_{j=1}^m$ is an antichain in a
finite poset $P$ (Figure~\ref{fig:poset}). An important and basic
fact, leading to HCF, was the observation \cite{BW83} that there
exists a nice and transparent formula describing the homotopy type
of the quotient $\Delta(P)/\Delta(P\setminus X)$ of order
complexes. Indeed it is elementary to see that,
\begin{equation}\label{eqn:difference}
\Delta(P)\setminus\Delta(P\setminus X) \cong \biguplus_{1\leq
j\leq m}{\rm OpenCone}(\Delta(P_{< x_j})\ast \Delta(P_{> x_j}))
\end{equation}
where the open cone ({\rm OpenCone}(Z)) with the base $Z$ is
defined as ${\rm Cone}(Z)\setminus Z$. By taking one-point
compactification of both sides of the homeomorphism
(\ref{eqn:difference}) we obtain the formula,
\begin{equation}\label{eqn:HCF-discrete}
\Delta(P)/\Delta(P\setminus X) \simeq \bigvee_{1\leq j\leq m}
\Sigma(\Delta(P_{< x_j})\ast \Delta(P_{> x_j})).
\end{equation}
Bj{\"o}rner and Walker in \cite{BW83} observed that if
$P\cup\{\widehat{0}, \widehat{1}\}$ ($P$ with added maximum and
minimum elements) is a lattice, and if the antichain $X = Co(y) :=
\{x\in P \mid x\vee y = \widehat{1}, x\wedge y = \widehat{0} \}$
arises as the set of all `complements' of a chosen element $y\in
P$, then the poset $P\setminus X$ is contractible. Then the
`Homotopy complementation formula' is the statement saying that
under these conditions,
\begin{equation}\label{eqn:HCF-discrete-final}
\Delta(P) \simeq \bigvee_{1\leq j\leq m} \Sigma(\Delta(P_{<
x_j})\ast \Delta(P_{> x_j})).
\end{equation}
When applied to the (truncated) lattice $\widetilde\Pi_n =
\Pi_n\setminus\{\widehat{0}, \widehat{1}\}$ of partitions of the
set $[n] = \{1,\ldots, n\}$ \cite{BW83}, the formula yields the
homotopy recurrence relation (\ref{eqn:hcf-1}) which immediately
leads to the computation of its homotopy type (described as a
wedge od spheres).

\begin{equation}\label{eqn:hcf-1}
\Delta(\widetilde\Pi_n) \simeq \bigvee_{i=2}^n
\Sigma(\Delta(\widetilde\Pi_{n-1}^i))
\end{equation}

\begin{equation}\label{eqn:hcf-2}
\Delta(\widetilde{\cal  G}_n(\mathbb{R})) \simeq S^{n-1}\wedge
\Sigma(\Delta(\widetilde{\cal  G}_{n-1}(\mathbb{R})))
\end{equation}

\begin{equation}\label{hcf-3}
\Delta(\widetilde{\cal  G}_n^{\pm}(\mathbb{\mathbb{R}})) \simeq
(S^{n-1}\vee S^{n-1})\wedge \Sigma(\Delta(\widetilde{\cal
G}^{\pm}_{n-1}(\mathbb{R})))
\end{equation}

\begin{equation}\label{eqn:hcf-4}
\Delta({\rm exp}_n(S^1)) \simeq S^{n}\wedge (\Delta({\cal B}_{n})/
\partial\Delta({\cal  B}_{n}))
\end{equation}

\begin{equation}\label{eqn:hcf-5}
\Delta({\cal  P}_n) \simeq P_n\wedge \Sigma(\Delta({\cal
P}_{n-1}))
\end{equation}

\begin{equation}\label{eqn:hcf-6}
\Delta({\rm exp}_n(X)) \simeq {\rm Thom}_{n}(X\setminus\{x_0\})
\end{equation}

\medskip\noindent
Recall that $\Pi_n$ is the lattice of all (unordered) partitions
of the set $[n] = \{1,\ldots, n\}$ (where $p_1 \prec p_2$ if $p_2$
is a refinement of $p_1$) and $\widetilde\Pi_n :=
\Pi_n\setminus\{\widehat{0}, \widehat{1}\}$.

\medskip The starting point of \cite{Zi98} was the observation
that similar ideas can be applied to the analysis of homotopy
types of order complexes of interesting topological posets. The
formulas (\ref{eqn:hcf-2}) to (\ref{eqn:hcf-6}) illustrating this
phenomenon are taken from \cite[Section~2]{Zi98}.

In order to establish a link from
(\ref{eqn:HCF-discrete})--(\ref{eqn:HCF-discrete-final}) to
(\ref{eqn:hcf-2})--(\ref{eqn:hcf-6}) let us take a look again at
Figure~\ref{fig:poset}. This time however we interpret $P$ as a
topological poset, so the antichain $X$ is a (not necessarily
discrete) topological space, while the decomposition
(\ref{eqn:difference}) of the space $\Delta(P)\setminus
\Delta(P\setminus X)$ is interpreted as a fibre bundle $\xi$ over
$X$.

\medskip
Moreover the space $\Delta(P)/\Delta(P\setminus X)$ is described
as a `Thom-space' (one-point compactification) of the bundle $\xi$
(see Proposition~4.8. and Corollaries~4.10.--4.12. in \cite{Zi98}
for more precise statements). If this bundle is trivial  the
Thom-space reduces to a smash product, as illustrated by the
schematic formula (\ref{eqn:hcf-5}), which subsumes both
(\ref{eqn:hcf-2}) and (\ref{hcf-3}). The relation
(\ref{eqn:hcf-2}) can be used for a proof of the homeomorphism
(\ref{eqn:Grass-poset}) (Remark~\ref{rem:remarkable}). The
relation (\ref{hcf-3}) provides a basis for a similar result about
the Grassmannian of oriented subspaces of $\mathbb{R}^n$.

\medskip
Suppose that $X$ is a finite $CW$-complex and let
$\mbox{exp}_n(X)$ be the topological poset of all non-empty
subsets of $X$ of size $\leq n$ (Example~3.3 in
\cite[Section~3]{Zi98}). If $x_0\in X$ is a chosen base point then
the set $Co(\{x_0\})$ of all complements of $\{x_0\}$ in
$\mbox{exp}_n(X)$ turns out to be the space $B(Y,n) = F(Y,n)/S_n$
of all unordered $n$-tuples in $Y := X\setminus\{x_0\}$. The
associated vector bundle $\xi$ is the canonical vector bundle,
\[
\mathbb{R}^{n-1}\longrightarrow F(Y,n)\times_{S_n}
V\longrightarrow B(Y,n)
\]
where $V\cong \mathbb{\mathbb{R}}^{n-1}$ is the standard
$(n-1)$-dimensional, permutation representation of the group
$S_n$. The associated `Thom-space' is the one-point
compactification
$$
{\rm Thom}_n(Y) := (F(Y,n)\times_{S_n} V)\cup \{\infty\} .
$$

The following result \cite[Theorem~5.8.]{Zi98} gives a complete
description of the homotopy type of the configuration poset
$\mbox{\rm exp}_n(X)$ in the category of admissible spaces
\cite[Definition~5.7.]{Zi98} (which include all finite
$CW$-complexes).

\begin{theo}\label{theo:rade-conf}
Suppose that $(X, x_0)$ is a finite $CW$-complex. Then,
\begin{equation}\label{eqn:thom-thom}
\Delta({\rm exp}_n(X)) \simeq \mbox{\rm
Thom}_n(X\setminus\{x_0\}).
\end{equation}
\end{theo}
\noindent The formula (\ref{eqn:thom-thom}) has a particularly
simple form if $X = S^1$ when $Y = X\setminus \{x_0\}\cong
\mathbb{R}^1$ and $B(Y,n)$ is homeomorphic to the interior of an
$n$-dimensional simplex. This has as a consequence the formula
(\ref{eqn:hcf-4}) (where ${\cal B}_{n}= \{I\subseteq [n] \mid
I\neq\emptyset\}$) which eventually leads to the proof that
$\Delta({\rm exp}_n(S^1))\cong S^{2n-1}$, see also \cite{Vas92,
Vas97} and \cite{Zi98} for more direct proofs.

\subsection{Weighted barycenter spaces and a conjecture of Vassiliev}
\label{sec:barycenter}

The following construction has been introduced by Vassiliev under
the name {\em simplicial  resolution} of configuration spaces.
Suppose that a smooth, compact manifold or more generally a finite
CW-complex $M$ is generically embedded in the space $\mathbb{R}^N$
of very large dimension $N$. Let ${\rm Conv}_r(M)$ be the union of
all (closed) $(r-1)$-dimensional simplices with vertices in the
embedded space $M$. The genericity of the embedding means that two
simplices $\mbox{conv}(A)$ and $\mbox{conv}(B)$ spanned by
different sets $A\neq B$ of vertices must have disjoint interiors.
The space ${\rm Conv}_r(X)$ is referred to as the {\em $r$-th
generic convex hull} of $X$.

\medskip
The following proposition records for the future reference a
simple fact that the order complex $\Delta({\rm exp}_n(M))$ can be
seen as the barycentric subdivision of the $n$-th generic convex
hull of $X$. Note that ${\rm Conv}_r(X)$ can be appropriately
described as a `continuous simplicial complex' on $X$ as the
(continuous) set of vertices.

\begin{prop}{\em (\cite[Section~5.2]{Zi98})}\label{prop:future}
Suppose that $X$ is a finite $CW$-complex. Then there is a natural
homeomorphism,
\begin{equation}\label{eqn:bary}
{\rm Conv}_n(X) \longrightarrow \Delta(\mbox{\rm exp}_n(X))
\end{equation}
of the $n$-th generic convex hull of $X$ and the order complex of
the corresponding configuration poset $\mbox{\rm exp}_n(X)$.
\end{prop}

\bigskip
The following conjecture, relating the order complex of $\mbox{\rm
exp}_n(X)$ and the $n$-fold, iterated join $X^{\ast n}$ of $X$,
was formulated by Victor Vassiliev on the conference ``Geometric
Combinatorics'' (MSRI Berkeley, February 1997).

\begin{equation}\label{eqn:vas-conjecture}
\Delta(\mbox{\rm exp}_n(X)) \simeq X\ast\ldots \ast X \cong
X^{\ast n}
\end{equation}
The conjecture (\ref{eqn:vas-conjecture}) was known to be true in
the case $X=S^1$ and this case  played a very important role in
applications.

\medskip
The whole `theory' of topological posets developed in \cite{Zi98}
was originally motivated by this conjecture. As a consequence of
Theorem~\ref{theo:rade-conf} it was shown
\cite[Proposition~5.10]{Zi98} that $\mbox{\rm exp}_n(S^2)$ does
not have the homotopy type of a sphere for $n\geq 2$, and in
particular,
$$
\Delta(\mbox{\rm exp}_n(S^2)) \neq (S^2)^{\ast n}
$$
which means that the conjecture is false already in the case of
the $2$-sphere.

\medskip
This result settled the general conjecture in the negative,
however this was not the end of the story. Kallel and Karoui,
motivated by some questions from non-linear analysis, began the
analysis of the space ${\rm Conv}_n(X)$ in \cite{KK11} from a
slightly different point of view. They used an alternative
description of this space as the space $\mathcal{B}_n(X)$ of all
{\em weighted barycenters} of $n$ or less points in $X$ ($=$ the
space of probability measures on $X$ with finite support of size
$\leq n$). Kallel and Karoui were familiar with the fact that this
space was used by Vassiliev, however they were apparently unaware
of \cite{Zi98}, in particular they were unaware of the original
Vassiliev's conjecture. Surprisingly enough, one of their main
results casts a new and interesting light on Vassiliev's
conjecture.

\begin{theo}\label{thm:Kallel-Karoui}{\em \cite[Theorem~1.1.]{KK11}}
Suppose that $X$ is a finite, connected $CW$-complex and let
$Sym^{\ast n}(X) := X^{\ast n}/S_n$ be the symmetric, $n$-fold
join of $X$. Then,
\begin{equation}\label{eqn:sym-join}
\mathcal{B}_n(X) \simeq Sym^{\ast n}(X).
\end{equation}
\end{theo}
In light of the homeomorphism $\Delta(\mbox{\rm exp}_n(X))\simeq
\mathcal{B}_n(X)$ (Proposition~\ref{prop:future}), and by
comparing (\ref{eqn:vas-conjecture}) and (\ref{eqn:sym-join}), we
see that after all Victor Vassiliev was right when he conjectured
that the homotopy type of the $n$-th generic convex hull of $X$ is
closely related to the iterated join of the space $X$.

\bigskip
Kallel and Karoui have established in \cite{KK11} many other
interesting results about spaces of weighted barycenters. For
example they establish a `symmetric smash product' formula for the
space  ${\mathcal B}_n(X)$.

\begin{theo}{\em (\cite[Theorem~5.3.]{KK11})}
$$
{\mathcal B}_n(X) \simeq S^{n-1} \wedge_{S_n} X^{(n)}.
$$
\end{theo}
As a consequence they deduced the following neat result of I.
James, E. Thomas, H. Toda, and J.H.C. Whitehead,
\begin{equation}
{\mathcal B}_2(S^n) \simeq \Sigma^{n+1}(\mathbb{R}P^n).
\end{equation}

It is not surprising that Theorem~\ref{theo:rade-conf} is equally
effective and elegant for computations of these examples. For
example $\Delta(\mbox{\rm exp}_2(S^n))$ has the same homotopy type
as the one-point compactification of
$$
F(\mathbb{R}^n,2)\times_{\mathbb{Z}/2} \mathbb{R}^1 \cong
\mathbb{R}^n \times
(\mathbb{R}^n\setminus\{0\})\times_{\mathbb{Z}/2} \mathbb{R}^1
\cong \mathbb{R}^n \times \mathbb{R}^+ \times
(S^{n-1}\times_{\mathbb{Z}/2} \mathbb{R}^1).
$$
Let $Z^+ = Z\cup\{\infty\}$ be the one-point compactification of a
locally compact space $Z$. Since $(U\times V)^+ \cong U^+ \wedge
V^+$ and $(S^{n-1}\times_{\mathbb{Z}/2} \mathbb{R}^1)^+ \cong
\mathbb{R}P^n$ we immediately observe that,
$$
{\mathcal B}_2(S^n) \cong \Delta(exp_2(S^n)) \simeq S^{n+1}\wedge
\mathbb{R}P^n \cong \Sigma^{n+1}(\mathbb{R}P^n).
$$
A similar argument based on Theorem~\ref{theo:rade-conf} can be
used for the proof of the homotopy equivalence,
$$
\Delta(exp_n(X)) \simeq S^{n-1} \wedge_{S_n} X^{(n)}.
$$

\section{Homotopy colimits and the index inequality}
\label{sec:sarkaria}

Perhaps the first appearance of homotopy colimits in combinatorial
applications was the application of this technique in \cite{ZZTop}
to the computation of (stable) homotopy types of arrangements of
subspaces, their links and complements. This paper was followed by
\cite{WZZ} and \cite{Zi98} and today diagrams of spaces and their
homotopy colimits are used more and more in geometric and
topological combinatorics, see \cite[Chapter 15]{K08} for a less
technical presentation directed towards combinatorially minded
readers.

\medskip
Formally a diagram of spaces over a finite poset $P$ is a functor
$\mathcal{D} : P \rightarrow Top$ from the poset category $P$ to
the category of topological spaces. Informally, a diagram over $P$
is a poset $P$ where each element $p\in P$ is associated a space
$D_p$ and for each pair $p\leq q$ there is a map $d_{pq}$
satisfying natural commutativity relations:
\[
\mbox{\rm For each} \, p\in P, \, d_{pp} = \mathbf{1}_{D_p} \,
  \mbox{\rm and for each triple} \, p\leq q\leq r,
\, d_{pq}\circ d_{qr} = d_{pr}.
\]
Each diagram can be associated a topological poset $P_\mathcal{D}$
where $P_\mathcal{D} = \coprod_{p\in P} D_p$ is the disjoint union
of all spaces $D_p$ (elements of $P_\mathcal{D}$ are pairs $(p,
x)$ where $x\in D_p$) and $(p,x)\prec (q,y)$ if and only if $p\leq
q$ and $d_{pq}(y) = x$. A nice consequence of this point of view
is the following relation,
\begin{equation}\label{eqn:diagrams}
\mathbf{hocolim}(\mathcal{D}) \cong \Delta(P_\mathcal{D})
\end{equation}
saying that the homotopy colimit in the case of diagram of spaces
over posets reduces essentially to the order complex construction
applied to topological posets.

\bigskip
The `Sarkaria's inequality', originally introduced and proved in
\cite{Ziv}, is one of the central results used in combinatorial
applications of equivariant index theory. The reader is referred
to \cite[Chapter 5]{Ma03} for a very nice exposition of this and
related results with numerous applications in topological
combinatorics. Recall that the index ${\rm Ind}_G(X)$ of a
$G$-space $X$ is a measure of complexity of $X$ which can be used
for proving Borsuk-Ulam type statements, for example the usual
Borsuk-Ulam theorem follows from the fact that ${\rm
Ind}_{\mathbb{Z}_2}(S^n) = n > {\rm
Ind}_{\mathbb{Z}_2}(S^{n-1})=n-1$.

\medskip
In general, for a given sequence $\mathcal{A} =
\{A_nG\}_{n=0}^{+\infty}$ of $G$-spaces, the associated
$\mathcal{A}$-index is defined by,
\begin{equation}\label{eqn:index}
Ind_G^{\mathcal{A}}(X) := \mbox{ {\rm Inf}}\{n\in \mathbb{N} \mid
X \stackrel{G}{\longrightarrow} A_nG\}
\end{equation}
where $X \stackrel{G}{\longrightarrow} Y$ means that there exists
a $G$-equivariant map from $X$ to $Y$.

\begin{prop}{\em (Sarkaria's inequality)}\label{prop:Sarkaria}
Let $G$ be a finite group and let $\mathcal{A} =
\{A_nG\}_{n=0}^{+\infty}$ be a sequence of $G$-spaces such that
$A_pG\ast A_qG \stackrel{G}{\longrightarrow}  A_{p+q+1}G$ for each
$p$ and $q$. Suppose that $L_0$ is a finite $G$-simplicial complex
and let $L\subset L_0$ be its $G$-invariant subcomplex. Then there
is an inequality,
\begin{equation}\label{eqn:Sark-inequality}
{\rm Ind}_G^{\mathcal{A}}(L) \geq {\rm Ind}_G^{\mathcal{A}} (L_0)
- {\rm Ind}_G^{\mathcal{A}} (\Delta(L_0\setminus L)) -1,
\end{equation}
where $\Delta(L_0\setminus L)$ is the order complex of the poset
$(L_0\setminus L, \subset)$.
\end{prop}
The proof of Proposition~\ref{prop:Sarkaria} is identical to the
proof given in \cite[Section 5.7]{Ma03} (and the original paper
\cite{Ziv}) so we leave the details to the interested reader.
\hfill $\square$

\medskip

Once the reader is prepared to emulate and extend the argument
used in the proof of Proposition~\ref{prop:Sarkaria} to the case
of topological posets the following proposition is a natural
consequence. We leave the details of the proof to the reader and
postpone the application of this extension of Sarkaria's
inequality to some other publication.

\begin{prop}\label{prop:Sarkaria-Zivaljevic}
Let $G$ be a finite group and let $\mathcal{A} =
\{A_nG\}_{n=0}^{+\infty}$ be a sequence of $G$-spaces such that
$A_pG\ast A_qG \stackrel{G}{\longrightarrow}  A_{p+q+1}G$ for each
$p$ and $q$. Suppose that $P$ is a finite (not necessarily free)
$G$-poset and let $P_0\subset P$ be its initial, $G$-invariant
subposet. Let $P_1 = P\setminus P_0$ be the complementary subposet
of $P$. Assume that $\mathcal{D} : P \rightarrow Top$ is a
$G$-diagram of spaces with $G$-action on $\mathcal{D}$ compatible
with the action on $P$ and let $\mathcal{D}_0$ and $\mathcal{D}_1$
be the restrictions of this diagram on $P_0$ and $P_1$
respectively. Then,
\begin{equation}\label{eqn:S-Ziv-inequality}
{\rm Ind}_G^{\mathcal{A}}(\| \mathcal{D}_0\|) \geq {\rm
Ind}_G^{\mathcal{A}}(\| \mathcal{D}\|) - {\rm
Ind}_G^{\mathcal{A}}(\| \mathcal{D}_1\|) -1,
\end{equation}
where $\| \mathcal{E}\| = \mathbf{hocolim}(\mathcal{E})$ is the
homotopy colimit of the diagram $\mathcal{E}$.
\end{prop}

\bigskip\noindent
{\bf Acknowledgements:} The referee's appropriate and thoughtful
remarks were of considerable help in improving the presentation of
results in the paper.


\small \baselineskip3pt
\begin{thebibliography}{10000}


\bibitem[AD12]{AD12} L. Anderson, E. Delucchi. Foundations for a theory
of complex matroids, arXiv:1005.3560v2 [math.CO].

\bibitem[BC88]{BC88} A. Bahri, J.M. Coron, On a non-linear elliptic equation
involving the critical sobolev exponent: the effect of the
topology of the domain, \textit{Comm. Pure and Applied
Mathematics}, Vol XLI, (1988), 253–-294.

\bibitem[BW83]{BW83}
A.~Bj{\"o}rner and J.W. Walker.
\newblock A homotopy complementation formula for partially ordered sets.
\newblock {\em European J. Combin.}, 4:11--19, 1983.

\bibitem[GKP]{GKP} R.L. Graham, D.E. Knuth, O. Patashnik.
\textit{Concrete Mathematics: A Foundation for Computer Science},
Second Edition, Addison-Wesley Professional 1994.

\bibitem[JV\v{Z}]{JVZ} D. Jojic, S. Vre\' cica, R. \v Zivaljevi\'
c. Symmetric multiple chessboard complexes and a new theorem of
Tverberg type,  arXiv:1502.05290v2 [math.CO].


\bibitem[KW08]{KW08} G. Kalai, A. Wigderson. Neighborly Embedded
Manifolds. \textit{Discrete and Computational Geometry}, vol. 40,
no. 3, pp. 319-324, 2008.

\bibitem[KK11]{KK11} S. Kallel, R. Karoui. Symmetric joins and weighted
barycenters, \textit{Advanced Nonlinear Studies} 11 (2011),
117--143. arXiv:math/0602283v3 [math.AT].


\bibitem[K08]{K08}  D. Kozlov. \textit{Combinatorial Algebraic
Topology}, Series `Algorithms and Computation in Mathematics',
Vol. 21, Springer 2008.

\bibitem [Ma03]{Ma03}
J.~Matou\v sek.
\newblock {\em Using the Borsuk-Ulam Theorem.}
\newblock {\em Lectures on Topological Methods in Combinatorics and Geometry.}
\newblock Springer-Verlag, Berlin, 2003.

\bibitem[KR97]{KR97} D.A. Klain, G-C. Rota. \textit{Introduction to Geometric
Probability}, Lezioni Lincee, Cambridge University Press 1997.

\bibitem[R98]{Rota98}
G-C. Rota.
\newblock Ten Mathematics Problems I will never solve.
\newblock {\em DMV mitteilungen} 2, 45--52, (1998).

\bibitem[S71]{S71} R. Schneider. On steiner points of convex
bodies, {\em Israel J.\ Math.}, 1971, Vol.~9, 241--249.

\bibitem[S79]{S79} R. Schneider. Boundary structure and curvature of convex
bodies; In Contributions to Geometry: Proceedings of the
Geometry-Symposium held in Singen June 28, 1979 to July 1, 1978/
eds. Ju\" urgen T\"{o}lke; J\"{o}rg M. Wills. Springer 1979.

\bibitem[R70]{R70} R. T. Rockafellar. {\em Convex Analysis}. Princeton Univ. Press
(1972).

\bibitem[Vas91]{Vas91}
V.A. Vassiliev.
\newblock Geometric realization of the homology of classical Lie
groups and complexes, S--dual to flag manifolds.
\newblock {\em St.--Petersburg Math. J.} 3:4, 108--115, (1991).

\bibitem[Vas92]{Vas92} V.A. Vassiliev. \textit{Complements of Discriminants of Smooth Maps:}
Topology and Applications: Revised Edition, A.M.S. 1992,
Translations of Mathematical Monographs, vol. 98.

\bibitem[Vas93]{Vas93} V.A. Vassiliev. Invariants of knots and
complements of discriminats. In \textit{``Developments in
Mathematics, the Moscow School''} (V.I.~Arnold, M.~Monasirsky,
eds.), Chapmann \& Hall, 1993, 194--250.

\bibitem[Vas97]{Vas97} V.A. Vassiliev. \textit{Topology of complements of discriminants}, Moscow, Phasis, 1997,
552 pp. (in Russian).


\bibitem[Vas99]{Vas99}
V.A. Vassiliev. Topological order complexes and resolutions of
discriminant sets. \textit{Publications de l'Institut
Math\'{e}matique}, (N.S.) 66 (80),  165--185, (1999).

\bibitem[WZ\v{Z}]{WZZ}
V.~Welker, G.M.~Ziegler, R.T.~\v Zivaljevi\' c.
\newblock Homotopy colimits--comparison lemmas for combinatorial
applications,
\newblock preprint 1997.


\bibitem[Z\v{Z}93]{ZZTop}
G.~M. Ziegler and R.~T.~{\v{Z}}ivaljevi{\'c}.
\newblock Homotopy types of subspace arrangements via diagrams of spaces.
\newblock {\em Math. Ann.}, 295:527--548, 1993.


\bibitem[Zie]{Zie} G.M. Ziegler. \textit{Lectures on Polytopes}.  Graduate Texts in Mathematics, Vol.
152, Springer-Verlag 1995.

\bibitem[\v{Z}i89]{Zi89}  R.T. \v Zivaljevi\' c. Extremal Minkowski additive selections of compact convex
sets. \textit{Proc. Amer. Math. Soc}, Vol. 105, (1989) 697--700.

\bibitem[\v{Z}-I-II]{Ziv}
R. \v Zivaljevi\' c.
\newblock User's guide to equivariant
methods in combinatorics, I and II.
\newblock {\it Publ. Inst. Math. (Beograd) (N.S.)}, (I) 59(73):114--130, 1996 and (II)
64(78):107--132, 1998.

\bibitem[\v{Z}i98]{Zi98}
R.T. \v Zivaljevi\' c.
\newblock Combinatorics of topological posets; Homotopy
complementation formulas.
\newblock {\em Advances in Applied mathematics} 21,
(1998) 547-–574.

\bibitem[\v{Z}i98b]{Zi98b} R.T. \v Zivaljevi\' c. Combinatorics of topological
posets. Lecture on the conference \textit{``Geometric
Combinatorics''}, Satellite conference of the International
Congress of Mathematics in Berlin 1998; Kotor, Yugoslavia, 28. 8.
-- 3. 9. 1998.
\url{http://www.ssag.matf.bg.ac.rs/konferencije/satellite/}.


\bibitem[\v{Z}i09]{Zi09} R.T. \v Zivaljevi\' c. Complex and quaternionic relatives of oriented matroids
(unpublished manuscript).
\end{thebibliography}
\end{document}